\documentclass{amsart}
\usepackage{amsmath}
\usepackage{amsfonts}
\usepackage{amssymb}
\usepackage{graphicx}%
\usepackage{mathrsfs}

\usepackage{times}
\usepackage{bbold}
\usepackage{stmaryrd}
\usepackage{array}
\usepackage{dsfont}
\usepackage[T1]{fontenc}

\newtheorem{theorem}{Theorem}

\newtheorem*{theoremMZ}{Theorem MZ}

\newtheorem*{mtheorem}{Theorem~\ref{T1} (Explicit)}
\newtheorem*{conjectureG-MZ}{Gaussian~vs.~MZ Derivatives Conjecture}
\newtheorem*{conjectureR-MZ}{Riemann~vs.~MZ Derivatives Conjecture}
\newtheorem*{SconjectureG-MZ}{Symmetric~Gaussian~vs.~MZ Derivatives Conjecture}
\newtheorem*{SconjectureR-MZ}{Symmetric~Riemann~vs.~MZ Derivatives Conjecture}

\newtheorem*{definitionMZ}{Definition MZ}
\newtheorem{example}[theorem]{Example}

\newtheorem*{problemMZ}{Problem MZ}
\newtheorem{proposition}[theorem]{Proposition}

\numberwithin{theorem}{section}

\begin{document}
\title[Counterexamples to the Gaussian vs. MZ Derivatives Conjecture]{Counterexamples to the Gaussian vs. MZ Derivatives Conjecture}

\author{J. Marshall Ash}
\address{Department of Mathematics, DePaul University\\Chicago, IL 60614}
\email{mash@depaul.edu}
\author{Stefan Catoiu}
\address{Department of Mathematics, DePaul University\\Chicago, IL 60614}
\email{scatoiu@depaul.edu}
\thanks{September 16, 2024. This paper is in final form and no version of it will be submitted for publication elsewhere.}
\subjclass[2010]{Primary 26A24; Secondary 05A30; 26A27; 47B39.}
\keywords{Gaussian derivative; generalized Riemann derivative; MZ derivative; Peano derivative; q-derivative; symmetric generalized Riemann derivative; symmetric Peano derivative.}

\begin{abstract}
J. Marcinkiewicz and A. Zygmund proved in 1936 that the special $n$-th generalized Riemann derivative ${_2}D_nf(x)$ with nodes $0,1,2,2^2,\ldots, 2^{n-1}$, is equivalent to the $n$-th Peano derivative $f_{(n)}(x)$, for all $n-1$ times Peano differentiable functions $f$ at~$x$. Call every $n$-th generalized Riemann derivative with this property an MZ derivative.
The recent paper Ash, Catoiu, and Fejzi\'c [Israel J. Math. {255} (2023):177--199] introduced the $n$-th Gaussian derivatives as the $n$-th generalized Riemann derivatives with nodes either $0,1,q,q^2,\ldots ,q^{n-1}$ or $1,q,q^2,\ldots ,q^{n}$, where~$q\neq0,\pm 1$, proved that the Gaussian derivatives are MZ derivatives, and conjectured that these are \emph{all} MZ derivatives.
In this article, we invalidate this conjecture by means of two counterexamples. The order in which these are presented allows an update of the conjecture after each counterexample. The proof of the first counterexample is simple, by scales of generalized Riemann derivatives. The proof of the second involves the classification of generalized Riemann derivatives of Ash, Catoiu, and Chin [Proc. Amer. Math. Soc {146} (2018):3847--3862].
Symmetric versions of all the results are also~included.
\end{abstract}
\maketitle

\noindent
\subsection{Peano and generalized Riemann differentiations}
For a positive integer $n$, an $n$-th \emph{generalized Riemann difference} of a function $f$ at $x$ is a difference of the form
\[
\Delta_{\mathcal{A}}(h,x;f)=\sum_{i=0}^ma_if(x+b_ih),
\]
where the datum $\mathcal{A}=\{a_i;b_i\}$ has non-zero $a_i$ and distinct $b_i$, and these satisfy the Vandermonde linear system $\sum_i a_i(b_i)^j=\delta_{j,n}\cdot n!$, for $j=0,1,\ldots ,n$. The numbers $x+b_ih$ are the \emph{base points}, and the numbers $b_i$ are the \emph{nodes} of the given difference. The above linear system is consistent when $m\geq n$, and has a unique solution when $m=n$, in which case the difference~is~\emph{exact}. The best known exact $n$-th generalized Riemann differences are the $n$-th \emph{Riemann difference},
\[
\Delta_n(h,x;f)=\sum_{i=0}^n(-1)^i\binom nif(x+(n-i)h),
\]
and the $n$-th \emph{symmetric Riemann difference},
\[
\Delta_n^s(h,x;f)=\sum_{i=0}^n(-1)^i\binom nif(x+(\frac n2-i)h).
\]
More exact $n$-th generalized Riemann differences include their $q$-analogues, the differences~${_q}\Delta_n$ with nodes $0,1,q,q^2\ldots ,q^{n-1}$ and ${_q}\bar\Delta_n$ with nodes $1,q,q^2\ldots ,q^{n}$, called \emph{Gaussian Riemann differences}, and ${_q}\Delta_n^s$ with nodes $0,\pm1,\pm q,\pm q^2,\ldots ,\pm q^{m-1}$ for $n=2m$ and $\pm1,\pm q,\pm q^2,\ldots ,\pm q^{m-1}$ for $n=2m+1$, called \emph{symmetric Gaussian Riemann difference}.
These were all computed in~\cite{ACF1}. Their expressions involve Gaussian binomial coefficients. In this paper we will simply call them \emph{Gaussian differences}.

An $n$-th \emph{generalized Riemann derivative} is a limit of a difference quotient of the form
\[
D_{\mathcal{A}}f(x)=\lim_{h\rightarrow 0}\Delta_{\mathcal{A}}(h,x;f)/h^n,
\]
where $\Delta_{\mathcal{A}}$ is an $n$-th generalized Riemann difference.
In particular, the $n$-th \emph{Riemann derivative} $D_nf(x)$ (resp. the $n$-th  \emph{symmetric Riemann derivative} $D_n^sf(x)$) correspond to~$\Delta_n$ (resp. $\Delta_n^s$), the $n$-th  \emph{Gaussian derivatives} ${_q}D_nf(x)$ and ${_q}\bar D_nf(x)$ (resp. the $n$-th  \emph{symmetric Gaussian derivative} ${_q}D_n^sf(x)$) correspond to ${_q}\Delta_nf(x)$ and~${_q}\bar \Delta_nf(x)$ (resp. to ${_q}\Delta_n^sf(x)$). When the derivative $D_{\mathcal{A}}f(x)$ exists, we say that $f$ is $n$ times \emph{generalized Riemann differentiable} (of type $\mathcal{A}$) at $x$, or \emph{$\mathcal{A}$-differentiable} at $x$.

A \emph{scale} by a non-zero $r$ of an $n$-th generalized Riemann difference $\Delta_{\mathcal{A}}(h,x;f)$ with data vector $\mathcal{A}=\{a_i;b_i\}$ is the $n$-th generalized Riemann difference
\[
\Delta_{\mathcal{A}_r}(h,x;f):=\tfrac 1{r^n}\Delta_{\mathcal{A}}(rh,x;f)
\]
with data vector $\mathcal{A}_r=\{a_i/r^n;rb_i\}$. For all $\mathcal{A}$ and $r$, a function $f$ is $\mathcal{A}$-differentiable at~$x$ if and only if $f$ is $\mathcal{A}_r$-differentiable at $x$ and $D_{\mathcal{A}_r}f(x)=D_{\mathcal{A}}f(x)$.

The more widely used notion of an \emph{$n$-th difference} is any non-zero scalar multiple of an~$n$-th generalized Riemann difference. Vice-versa, an $n$-th generalized Riemann difference is a normalized $n$-th difference, where the normalization condition, or the $n$-th Vandermonde relation $\sum a_i(b_i)^n=n!$, forces the value of the generalized derivative to agree with the value of the ordinary $n$-th derivative, for all $n$ times ordinary differentiable functions at $x$.

An $n$th difference $\Delta_{\mathcal{A}}(h,x;f)$ is \emph{symmetric}, if $\Delta_{\mathcal{A}}(-h,x;f)=(-1)^n\Delta_{\mathcal{A}}(h,x;f)$. For example, both $\Delta_n^s$ and~${_q}\Delta_{n}^s$ are symmetric for all $n$. A generalized Riemann derivative~$D_{\mathcal{A}}$ is symmetric, if its difference~$\Delta_{\mathcal{A}}$ is symmetric.

Another family of generalized derivatives are the Peano derivatives. A function $f$ is $n$ times \emph{Peano differentiable} at $x$ if it is approximated to order $n$ about $x$ by a polynomial, i.e., if there are constants $c_0=f(x),c_1,\ldots ,c_n$ such that
\[
f(x+h)=c_0+\frac {c_1}{1!}h+\frac {c_2}{2!}h^2+\cdots +\frac {c_n}{n!}h^n+o(h^n).
\]
The numbers $c_0,c_1,\ldots ,c_n$ are usually denoted by $f_{(0)}(x),f_{(1)}(x),\ldots ,f_{(n)}(x)$ and called the first~$n$ \emph{Peano derivatives} of~$f$ at~$x$. Note that the existence of the $n$-th Peano derivative $f_{(n)}(x)$ implies the existence of all lower order Peano derivatives at $x$. Moreover, $f_{(0)}(x)=f(x)$ and $f_{(1)}(x)=f'(x)$, while, in higher orders, the $n$-th Peano derivative is known to be strictly more general than the ordinary derivative.

A function $f$ is $n$ times \emph{symmetric Peano differentiable} at $x$ if there exist real constants~$c_0,c_1,\ldots ,c_n$ such that
\[
\tfrac 12[f(x+h)+(-1)^nf(x-h)]=c_0+\frac {c_1}{1!}h+\frac {c_2}{2!}h^2+\cdots +\frac {c_n}{n!}h^n+o(h^n),
\]
in which case the numbers $c_i$, for $i=0,1,\ldots ,n$, are are denoted by $f_{(i)}^s(x)$ and called the first~$n$ \emph{symmetric Peano derivatives} of~$f$ at~$x$. By replacing $h$ with $-h$, we deduce that~$f_{(n-1)}^s(x)=f_{(n-3)}^s(x)=\cdots =f_{(1\text{ or }0)}^s(x)=0$, and that the existence of the~$n$th symmetric Peano derivative $f_{(n)}^s(x)$ implies the existence of all same-parity-lower-order~symmetric Peano derivatives of $f$ at $x$. The condition $c_0:=f(x)$ when $n$ is even is added to the above definition, in order to ensure that the existence of $f_{(0)}^s(x)$ is the same as the continuity of $f$ at $x$.

Two generalized differentiations are said to \emph{imply} (or to be \emph{equivalent} to) each other, if the existence of one implies (or is equivalent to) the existence of the other, for all functions~$f$ at~$x$. By Taylor or Peano expansion of the generalized Riemann difference, one can show that every $n$ times Peano differentiable function $f$ at~$x$ is~$n$ times generalized Riemann differentiable at~$x$ for any (admissible)~$\mathcal{A}$ and~$D_{\mathcal{A}}f(x)=f_{(n)}(x)$. The converse of this, or the equivalence between Peano and (generalized) Riemann differentiations, has been a problem since 1927, initiated by Khintchine in \cite{Ki}. Marcinkiewicz and Zygmund have shown in \cite{MZ} that the $n$-th symmetric Riemann derivative $D_n^s$ is equivalent to the $n$-th Peano derivative, for all functions $f$ at a.e. points $x$ on a measurable set, and Ash has shown in \cite{As} that the same result holds true for all $n$-th generalized Riemann derivatives. The measurability condition on the set was removed by Fejzi\'c and Weil in \cite{FW1}.

\subsection{MZ derivatives and two conjectures related to them}\label{S0.1}
Our first motivation comes from the following theorem, proved by Marcinkiewicz and Zygmund in 1936, that both of the above almost everywhere results relied on:

\begin{theoremMZ}[\cite{MZ}, $\S 10$, Lemma 1]
Let ${_2}D_n$, $n\geq 1$, be the $n$-th generalized Riemann derivative with nodes $0,1,2,4,\ldots ,2^{n-1}$. Then, for all functions $f$ and points $x$,
\[
\text{both $f_{(n-1)}(x)$ and ${_2}D_nf(x)$ exist $\Longleftrightarrow $ $f_{(n)}(x)$ exists.}
\]
\end{theoremMZ}

\noindent
In other words, the $n$-th generalized Riemann derivative ${_2}D_nf(x)$ is equivalent to the $n$-th Peano derivative $f_{(n)}(x)$, for all $n-1$ times Peano differentiable functions $f$ at $x$.

The present article is about $n$-th generalized Riemann derivatives $D_{\mathcal{A}}f(x)$ that can play the role of ${_2}D_nf(x)$ in Theorem~MZ. In this way, Theorem MZ prompts the following definition.



\begin{definitionMZ}
{\rm
An $n$-th generalized Riemann differentiation $D_{\mathcal{A}}$ is an \emph{MZ differentiation} at a point $x$, if, for all $n-1$ times Peano differentiable functions $f$ at $x$,
\begin{itemize}
\item[] $f$ is $n$ times Peano differentiable at~$x$ $\Longleftrightarrow $ $f$ is $\mathcal{A}$-differentiable~at~$x$.
\end{itemize}
}
\end{definitionMZ}

One basic problem of MZ differentiation is the following:
\begin{problemMZ}
{\rm
For each order $n$ and point $x$, determine all $n$-th MZ differentiations at $x$.}
\end{problemMZ}

A solution to Problem MZ was proposed in \cite{ACF1}, where
it was proved that, for all real numbers $q$, with $q\neq 0,\pm 1$, the $n$-th Gaussian derivatives ${_q}D_nf(x)$ with nodes $0,1,q,\ldots ,q^{n-1}$ and ${_q}\bar D_nf(x)$ with nodes $1,q,\ldots ,q^{n}$ are MZ derivatives at $x$. Furthermore, it was conjectured that these are all MZ derivatives of order~$n$ at $x$, that is,

\begin{conjectureG-MZ}[\cite{ACF1}, Conjecture A]{\rm Each MZ differentiation is a Gaussian differentiation.}
\end{conjectureG-MZ}

In other words, the Gaussian vs. MZ Derivatives Conjecture asserts that all MZ derivatives classify as Gaussian.
Article \cite{ACF1} proved that the conjecture is true for $n=1,2$, and left it open for~$n\geq 3$.

A particular case of the Gaussian vs. MZ Derivatives Conjecture is the following conjecture from \cite{AC2} on the relation between Riemann differentiation and MZ differentiation:

\begin{conjectureR-MZ}[\cite{AC2}, Conjecture~4.2] \,{\rm
\begin{itemize}
\item[] When $n\geq 3$, $D_n$ is not an MZ differentiation, for all functions $f$ at $x$.
\end{itemize}
}
\end{conjectureR-MZ}

This conjecture was first formally stated in \cite{AC2}. It was proved for $n=3$ in \cite{ACF} and for~$n=7$ in \cite{ACF1}. The Gaussian vs. MZ Derivatives Conjecture was proved very recently for general $n$ in \cite{CF1} and is now a theorem.\vspace{-.05in}
\[
\ast\quad\ast\quad\ast\vspace{-.05in}
\]
The goal of this article is to invalidate the Gaussian vs. MZ Derivatives Conjecture, and we do this by means of two counterexamples. The order we chose to arrange them allows a restatement of the conjecture in sharper form after each counterexample.

The first counterexample, given in Example~\ref{E1}, is a scale of a Gaussian derivative that is not a Gaussian. The conjecture is then updated to \emph{``every MZ derivative is a scale of a Gaussian''}. Proving that the derivative is not a Gaussian is easily done by inspection.~The motivation for this counterexample comes from \cite{C}, where all exact first order~MZ derivatives are determined, and all such derivatives are scales of Gaussian derivatives, but not all of them are Gaussian. The section ends with Proposition~\ref{P8}, characterizing when two~Gaussian derivatives are scales of each other.

The second counterexample, given in Example~\ref{E2}, is an MZ derivative equivalent to a Gaussian, but not a scale of a Gaussian. The conjecture is then further updated to \emph{``every MZ derivative is equivalent to a Gaussian''}. Proving that the derivative is not a scale of a Gaussian is also easily done by inspection. Proving that it is not equivalent to a Gaussian is done by using a powerful result, the classification of generalized Riemann derivatives of~\cite[Theorem~2]{ACCh}, explained in Section~\ref{S2.1} and restated as Theorem~\ref{T1}. Proposition~\ref{P3} asserts that, for exact generalized Riemann derivatives, the second update of the Gaussian vs. MZ derivatives conjecture is no different than the first one. This is a consequence of Theorem~\ref{T2.4}, characterizing all exact generalized Riemann derivatives equivalent to Gaussian derivatives.

Most of the results in Sections~1-2 have symmetric analogues. These are discussed in Section~\ref{S6} and the similar proofs are either sketched or omitted.\vspace{-.05in}
\[
\ast\quad\ast\quad\ast\vspace{-.05in}
\]
The Riemann derivatives were introduced by Riemann in the mid 1800s (see \cite{R}). The Peano derivatives were invented by Peano in \cite{P} in 1892 and developed by De la Vall\'ee Poussin in \cite{dlVP}. The generalized Riemann derivatives were introduced by Denjoy in \cite{D} in 1935. The symmetric generalized Riemann derivatives and the symmetric Peano derivatives that we will define in Section~\ref{S6} were recently investigated in \cite{AC2}. Other $q$-analogues of the Riemann derivatives were studied in \cite{AC,ACR}.

The generalized Riemann derivatives were shown to satisfy properties similar to the ones satisfied by the ordinary derivatives, such as monotonicity \cite{HL1,T,W}, convexity \cite{GGR1,HL,MM}, or the Mean Value Theorem \cite{AJ,FFR}. They have many applications in the theory of trigonometric series \cite{SZ,Z} and in numerical analysis \cite{AJJ,L,Sa}. Multidimensional Riemann derivatives are studied in \cite{AC1}. For more on generalized Riemann derivatives, see the expository article \cite{As2} by Ash.
The Peano derivatives also have a long and rich history; see the expository article \cite{EW} by Evans and Weil. More recent developments on Peano derivatives are found in \cite{F, F1,FR,FW,LPW}.

\section{The first counterexample.}
The first counterexample to the Gaussian vs. MZ Derivatives Conjecture provides a scale of a Gaussian derivative that is not a Gaussian derivative.

\begin{example}\label{E1}
{\rm
For all real $q$, $q\neq 0,\pm 1,\pm 1/\sqrt{2}$, the scale by 2 of the second Gaussian difference\vspace{-.05in}
\[{_q}\bar\Delta_2(h,x;f)=\tfrac {2!}{(q^2-1)(q^2-q)}[f(x+q^2h)-(q+1)f(x+qh)+qf(x+h)]\]
is the non-Gaussian second difference\vspace{-.05in}
\[\tfrac 14\cdot \tfrac {2!}{(q^2-1)(q^2-q)}[f(x+2q^2h)-(q+1)f(x+2qh)+qf(x+2h)].\]

Indeed, if this was a Gaussian difference, then one node would be 1, leading to $2q^2=1$ or~$2q=1$, both with impossible solutions.
}
\end{example}

A related result is the following proposition, describing all pairs of Gaussian differences that are either scalar multiples or scales of each other.

\begin{proposition}\label{P8} {\rm \,(i)} If two distinct Gaussian differences are scalar multiples of each other, then they must coincide.

{\rm \,(ii)} Two distinct Gaussian differences are scales of each other if and only if they have the form\vspace{-.1in}
\[\text{${_q}\Delta_n$ and either one of ${_{-q}}\Delta_n$ and ${_{\pm q^{-1}}}\Delta_n$,\quad for $n\geq 2$,}\vspace{-.1in}\]
or\vspace{-.1in}
\[\text{${_q}\bar\Delta_n$ and either one of ${_{-q}}\bar\Delta_n$ and ${_{\pm q^{-1}}}\bar\Delta_n$,\quad for $n\geq 1$.}\]
\end{proposition}

\begin{proof} (i) This follows from the same property for generalized Riemann derivatives, which is a consequence of the Vandermonde relations.

(ii) Suppose two Gaussian differences are scales of each other. Then either 0 is a node in both or is a node in neither. In the first case, the two differences are of the form ${_q}\Delta_n$ and~${_Q}\Delta_n$, for some $n$, $q$ and $Q$, with $q,Q\neq 0,\pm 1$ and $q\neq Q$. Their respective sets of nodes are $\{0,1,q,\ldots,q^{n-1}\}$ and $\{0,1,Q,\ldots,Q^{n-1}\}$. Let a scale by~$r$ map the first difference into the other. This maps the first set of nodes, in order, into~$\{0,r,rq,\ldots, rq^{n-1}\}$. The absolute values of the non-zero elements of the same (differently ordered) set are two monotonic sequences, $1,|Q|,\ldots,|Q^{n-1}|$ and~$|r|,|rq|,\ldots ,|rq^{n-1}|$. This forces either~$|r|=1$ and~$|rq|=|Q|$, or $|r|=|Q^{n-1}|$ and~$|rq|=|Q^{n-2}|$, that is, either~$Q=-q$ or~$Q=\pm q^{-1}$. The second case has a similar proof.
\end{proof}

Both Example~\ref{E1} and Proposition~\ref{P8} show that, in general, scales of Gaussian derivatives are not Gaussian, while, as we have seen earlier in the introduction, they are equivalent to Gaussian and so they are MZ derivatives. In particular, Example~\ref{E1} is a counterexample to the Gaussian vs. MZ Derivatives Conjecture.

The first update of the Gaussian conjecture, for which Example~\ref{E1} is not a counterexample, would either have to extend the notion of a Gaussian derivative to include all of its scales, or to reassert the same statement of the conjecture up to a scale.
We chose the second option as the first update of the Gaussian vs. MZ Derivatives Conjecture:

\begin{conjectureG-MZ}[The first update]\label{Cj2}
{\rm Each MZ differentiation is~a scale of a Gaussian differentiation.}
\end{conjectureG-MZ}

\section{The second counterexample}
The second counterexample to the Gaussian vs. MZ Derivatives Conjecture is investigated~in~Section~\ref{S2.2} using the classification of generalized Riemann derivatives explained in Section~\ref{S2.1}.
\subsection{The Classification of Generalized Riemann Derivatives}\label{S2.1}
Before getting into the details of this subsection, consider the following three first-order generalized Riemann derivatives:\vspace{-.05in}
{\small\[
D_{\mathcal{A}}f(x)=\lim_{h\rightarrow 0} \frac {f(x+h)-f(x-h)}{2h},\quad D_{\mathcal{B}}f(x)=\lim_{h\rightarrow 0} \frac {f(x+h)-f(x)}h,
\]}and\vspace{-.1in}
{\small\[
D_{\mathcal{C}}f(x)=\lim_{h\rightarrow 0} \frac {3f(x+h)-5f(x)+2f(x-h)}h.
\]}The first is the symmetric derivative~$f'_s(x)$, and the second is the ordinary derivative $f'(x)$. The example of $f(x)=|x|$, for which $f'_s(0)$ exists and is zero and $f'(0)$ does not exist, shows that $D_{\mathcal{A}}$ is not equivalent to $D_{\mathcal{B}}$ at $x=0$. And by taking the limit as $h\rightarrow 0$ in the identity\vspace{-.05in}
{\small\[
\frac 35\cdot\frac {3f(h)-5f(0)+2f(-h)}h+\frac 25\cdot\frac {3f(-h)-5f(0)+2f(h)}{-h}=\frac {f(h)-f(0)}h,
\]}the derivative $D_{\mathcal{C}}$ implies, hence is equivalent to, the derivative $D_{\mathcal{B}}$ at $x=0$. We deduce that in proving or disproving that two given generalized Riemann derivatives are equivalent one needs either some inventive algebra or a counterexample.

\medskip
Our second motivation comes from The Classification of Generalized Riemann Derivatives, due to Ash, Catoiu, and~Chin in~\cite{ACCh}, which was extended to complex functions in \cite{ACCH}. It characterizes all pairs of generalized Riemann differences whose associated generalized Riemann derivatives are equivalent for all functions $f$ at $x$.

A difference $\Delta_{\mathcal{A}}(h,x;f)$ is \emph{even} or \emph{odd}, if $\Delta_{\mathcal{A}}(-h,x;f)=\pm \Delta_{\mathcal{A}}(h,x;f)$.
The \emph{symmetrizer} and the \emph{skew-symmetrizer} of an $n$-th difference $\Delta_{\mathcal{A}}(h,x;f)$ are the differences
\[
\Delta_{\mathcal{A}}^{\pm }(h,x;f)=\tfrac 12[\Delta_{\mathcal{A}}(h,x;f)\pm (-1)^n\Delta_{\mathcal{A}}(-h,x;f)].
\]
The first is an $n$-th symmetric generalized Riemann difference, while the second is either zero or a symmetric difference whose order is higher than $n$ and has a different parity than~$n$. The identity $\Delta_{\mathcal{A}}(h,x;f)=\Delta_{\mathcal{A}}^+(h,x;f)+\Delta_{\mathcal{A}}^-(h,x;f)$ is the unique decomposition of $\Delta_{\mathcal{A}}(h,x;f)$ as a sum of even and odd differences.
A difference is symmetric if and only if it is its own symmetrizer, or if and only if its skew-symmetrizer is zero.

\medskip
Here is an equivalent statement of the result. 

\begin{theorem}[\cite{ACCh}, Theorem~2]\label{T1}
Let~$\mathcal{A}$~and~$\mathcal{B}$ be the data sets for two generalized Riemann derivatives of orders $n$ and $N$. Then, for a fixed real number $x$, the following assertions are equivalent:
\begin{enumerate}
\item[(i)] $D_{\mathcal{A}}f(x)$ is equivalent to $D_{\mathcal{B}}f(x)$ for every real function $f$;
\item[(ii)] $N=n$ and there exist non-zero constants $r,s,B$ such that
\[
\qquad\text{$\Delta_{\mathcal{A}}^{+}(h,x;f)=r^{-n}\cdot \Delta_{\mathcal{B}}^{+}(rh,x;f)$
and $\Delta_{\mathcal{A}}^{-}(h,x;f)=B\cdot\Delta_{\mathcal{B}}^{-}(sh,x;f)$,}
\]
for all $f$ and $h$.
\end{enumerate}
\end{theorem}
The proof of Theorem~\ref{T1} given in \cite{ACCh} uses an unusual tool from abstract algebra: the theory of ideals in a group algebra. Other instances where the theory of ideals is involved~in functional analysis are found in \cite{JPS,SZ1}.

Going back to the above example and with the notation $\Delta(h)$ for a difference $\Delta(h,0;f)$, the symmetrizers and skew-symmetrizers of the differences corresponding to the given derivatives are
{\small\[
\begin{array}{lll}
\Delta_{\mathcal{A}}^+(h)=\Delta_{\mathcal{A}}(h), & \Delta_{\mathcal{B}}^+(h)=\frac 12[f(h)-f(-h)], & \Delta_{\mathcal{B}}^-(h)=\frac 12[f(h)+f(-h)-2f(0)],\\
\Delta_{\mathcal{A}}^-(h)=0, & \Delta_{\mathcal{C}}^+(h)=\frac 12[f(h)-f(-h)], & \Delta_{\mathcal{C}}^-(h)=\frac 52[f(h)+f(-h)-2f(0)].
\end{array}
\]}By Theorem~\ref{T1}, $D_{\mathcal{A}}$ is not equivalent to $D_{\mathcal{B}}$ at $x=0$, since $\Delta_{\mathcal{A}}^-(h)=0$ and $\Delta_{\mathcal{B}}^-(h)\neq 0$; and $D_{\mathcal{C}}$ is equivalent to $D_{\mathcal{B}}$ at $x=0$, since $\Delta_{\mathcal{C}}^+(h)=\Delta_{\mathcal{B}}^+(h)$ and $\Delta_{\mathcal{C}}^-(h)=5\Delta_{\mathcal{B}}^-(h)$. Thus Theorem~\ref{T1} reduced the problem of deciding whether two generalized derivatives are equivalent to the problem of inspecting the symmetrizers and skew-symmetrizers of the two associated differences.
\subsection{The second counterexample}\label{S2.2}
The second counterexample to the Gaussian vs. MZ Derivatives Conjecture is a 
counterexample to the first update of the conjecture. This is based on using Theorem~\ref{T1} to construct a generalized Riemann derivative that is both equivalent to a Gaussian derivative and not a scale of a Gaussian derivative.
For simplicity, we denote each difference $\Delta (h,0;f)$ by $\Delta (h)$.

\begin{example}\label{E2}
{\rm
Let $q$ be a real number, with $q\neq 0,\pm 1$. Then the second Gaussian difference\vspace{-.05in}
\[
\Delta_{\mathcal{A}}(h):={_q}\bar\Delta_2(h)=\bar\lambda_2[f(q^2h)-(q+1)f(qh)+qf(h)],
\]
where $\bar\lambda_2={2!}/({(q^2-1)(q^2-q)})$, has the symmetrizer and skew-symmetrizer given by
{\small
\[
\begin{aligned}
\Delta_{\mathcal{A}}^{+}(h)&=\tfrac 12\bar\lambda_2[f(q^2h)-(q+1)f(qh)+qf(h)+qf(-h)-(q+1)f(-qh)+f(-q^2h)],\\
\Delta_{\mathcal{A}}^{-}(h)&=\tfrac 12\bar\lambda_2[f(q^2h)-(q+1)f(qh)+qf(h)-qf(-h)+(q+1)f(-qh)-f(-q^2h)].
\end{aligned}
\]}By taking $\Delta_{\mathcal{B}}(h):=\Delta_{\mathcal{A}}^{+}(h)+3\Delta_{\mathcal{A}}^{-}(h)$, or
{\small
\[
\Delta_{\mathcal{B}}(h)=\bar\lambda_2[2f(q^2h)-2(q+1)f(qh)+2qf(h)-qf(-h)+(q+1)f(-qh)-f(-q^2h)],
\]}by Theorem~\ref{T1}, the derivative $D_{\mathcal{A}}f(x)$ is equivalent to $D_{\mathcal{B}}f(x)$, while, by inspection, the second order generalized Riemann difference $\Delta_{\mathcal{B}}(h)$ is not a scale of a Gaussian difference.
}
\end{example}

At this point, we can further update the Gaussian conjecture by either extending the notion of a Gaussian derivative to contain not only all of its scales, but also all generalized Riemann derivatives that are equivalent to it; or, by reasserting the conjecture as its previous statement up to the equivalence of generalized Riemann derivatives.

We chose the second option as the following higher update of the Gaussian conjecture:

\begin{conjectureG-MZ}[The second update]\label{Cj3} {\rm Each MZ differentiation is equivalent to a Gaussian differentiation.}
\end{conjectureG-MZ}

As the above example is a non-exact generalized Riemann derivative, one might think that by restricting the definition of an MZ derivative to only exact generalized Riemann derivatives, the second update of the Gaussian vs. MZ Derivatives Conjecture is no different than the first update.
The next proposition shows that this is, indeed, the case.

\begin{proposition}\label{P3}
The second update of the Gaussian vs. MZ Derivatives Conjecture is the same as the first update, when all MZ derivatives are assumed exact.
\end{proposition}

This proposition is an immediate consequence of the following theorem, characterizing all exact generalized Riemann derivatives equivalent to a Gaussian.

\begin{theorem}\label{T2.4}
An exact generalized Riemann derivative is equivalent to a Gaussian if and only if it is a scale of a Gaussian.
\end{theorem}

\begin{proof}
Suppose $D_{\mathcal{B}}={_q}\bar D_n$ and let $D_{\mathcal{A}}$ be an exact $n$th generalized Riemann derivative equivalent to it. Then $\Delta_{\mathcal{A}}^{\pm}$ and $\Delta_{\mathcal{B}}^{\pm}$ satisfy the system of equations in Theorem~\ref{T1}(ii). Using a scale by $1/r$ of $\Delta_{\mathcal{B}}$ and without loss of generality, we may assume that $r=1$, hence the system is
\begin{equation}\label{eq1}
\qquad\text{$\Delta_{\mathcal{A}}^{+}(h,x;f)=\Delta_{\mathcal{B}}^{+}(h,x;f)$
,\quad $\Delta_{\mathcal{A}}^{-}(h,x;f)=B\cdot\Delta_{\mathcal{B}}^{-}(sh,x;f)$,}
\end{equation}
for some non-zero $B$ and $s$.

It suffices to show that $B,s=\pm 1$. Indeed, since the nodes of $\Delta_{\mathcal{B}}$ are~$1,q,\ldots,q^n$, those of $\Delta_{\mathcal{B}}^{\pm}$ and $\Delta_{\mathcal{A}}^{+}$ are $\pm1,\pm q,\ldots,\pm q^n$. In particular, at least one out of each pair~$\pm q^i$, for~$i=0,1,\ldots ,n$, is a node in $\Delta_{\mathcal{A}}$. This and $\Delta_{\mathcal{A}}$ being an exact $n$th generalized Riemann difference imply that the set of nodes of $\Delta_{\mathcal{A}}$ consists precisely of one node out of each pair~$\pm q^i$, for~$i=0,1,\ldots ,n$. The equality of the sets of nodes of the differences on both sides of the second equation in \eqref{eq1} yields $s=\pm 1$.

Assume $s=1$ and prove that either $B=1$ and $\Delta_{\mathcal{A}}=\Delta_{\mathcal{B}}$, or $B=-1$, $n$ is even, and $\Delta_{\mathcal{A}}$ is a scale by $-1$ of $\Delta_{\mathcal{B}}$. To see this, we start by writing the difference $\Delta_{\mathcal{B}}$ explicitly~as\vspace{-.05in}
\[
\Delta_{\mathcal{B}}(h,x;f)=\sum_{i=0}^nB_if(x+q^ih).
\]
Then, by \eqref{eq1} with $s=1$, the difference $\Delta_{\mathcal{A}}=\Delta_{\mathcal{A}}^++\Delta_{\mathcal{A}}^-=\Delta_{\mathcal{B}}^++B\Delta_{\mathcal{B}}^-$ has the explicit form\vspace{-.05in}
\[
\begin{aligned}
\Delta_{\mathcal{A}}&=\sum_{i=0}^nB_i\frac {f(x+q^ih)+f(x-q^ih)}2+B\sum_{i=0}^nB_i\frac {f(x+q^ih)-f(x-q^ih)}2\\
&=\frac 12\sum_{i=0}^nB_i\left\{(1+B)f(x+q^ih)+(1-B)f(x-q^ih)\right\}.
\end{aligned}
\]
Since, for each $i$, $\Delta_{\mathcal{A}}$ has precisely one node of the form $\pm q^{i}$ and $B_i$ is nonzero, we deduce that $1\pm B=0$, or $B=\pm 1$. The case $B=-1$ leads to $\Delta_{\mathcal{A}}(h,x;f)=\Delta_{\mathcal{B}}(-h,x;f)$, which is a generalized Riemann difference precisely when $n$ is even, that is, when $\Delta_{\mathcal{A}}$ is the scale by $-1$ of $\Delta_{\mathcal{B}}$. The case $B=1$ leads to~$\Delta_{\mathcal{A}}=\Delta_{\mathcal{B}}$. The case $s=-1$ is similar.
The proof for $D_{\mathcal{B}}={_q}D_n$ is similar to the proof for $D_{\mathcal{B}}={_q}\bar D_n$.
\end{proof}
\section{Symmetric Gaussian and MZ Derivatives}\label{S6}
\subsection{Symmetric MZ derivatives and two conjectures about them}
The symmetric version of Theorem~MZ is the following theorem proved in \cite{AC2}. For~$n$ at least two, we denote $m=\lfloor (n-1)/2\rfloor $ and define ${_2}D_n^sf(x)$ to be the $n$-th generalized~Riemann derivative with nodes $(0),\pm 1,\pm 2, \pm 4,\ldots ,\pm 2^{m}$, where $(0)$ means that 0 is taken only for~$n$ even.

\begin{theoremMZ}[Symmetric] \cite[Theorem~2.2]{AC2} For $n\geq 2$, let ${_2}D_n^sf(x)$ be the above $n$-th symmetric generalized Riemann derivative.
Then, for all functions $f$ and points $x$,
\[
\text{both $f_{(n-2)}^s(x)$ and ${_2}D_n^sf(x)$ exist $\Longleftrightarrow $ $f_{(n)}^s(x)$ exists.}
\]
\end{theoremMZ}

\medskip
This theorem leads to the following definition and problem:

\begin{definitionMZ}[Symmetric]
{\rm
An $n$-th symmetric generalized Riemann differentiation~$D_{\mathcal{A}}$ is a \emph{symmetric MZ differentiation} at $x$, if, for all functions $f$, 
\[
\text{both $f_{(n-2)}^s(x)$ and ${D}_{\mathcal{A}}f(x)$ exist $\Longleftrightarrow $ $f_{(n)}^s(x)$ exists.}
\]
}
\end{definitionMZ}
The following is the symmetric version of Problem~MZ:
\begin{problemMZ}[Symmetric]
{\rm
For each order $n$ and point $x$, determine all $n$-th symmetric~MZ differentiations at $x$.}
\end{problemMZ}

A solution to the symmetric Problem MZ was also proposed in \cite{ACF1}, where
it was proved that, for all real numbers $q$, with $q\neq 0,\pm 1$, the $n$-th symmetric Gaussian derivatives~${_q}D_n^sf(x)$ are symmetric MZ derivatives at $x$, and then it was conjectured that these are all symmetric MZ derivatives at $x$ of order~$n$, that is,

\begin{SconjectureG-MZ}[\cite{ACF1}, Conjecture A]{\rm In orders at least 3, each symmetric MZ differentiation is a symmetric Gaussian differentiation.}
\end{SconjectureG-MZ}

\noindent
In other words, all symmetric MZ derivatives classify as symmetric Gaussian.
Article~\cite{ACF1} proved it to be false for $n=1,2$, true for $n=3,4$, and left it open for~$n\geq 5$.

A particular case of the symmetric Gaussian vs. MZ Derivatives Conjecture is the following conjecture from \cite{AC2} on the relation between symmetric Riemann and symmetric~MZ differentiations:

\begin{SconjectureR-MZ}[\cite{AC2}, Conjecture~4.1]\,{\rm For $n\geq 5$, the $n$-th symmetric Riemann differentiation $D_n^s$ is not a symmetric MZ differentiation.}
\end{SconjectureR-MZ}

This conjecture was first stated in \cite{AC2} for $n\geq 3$. It is obviously false for $n=2$. Article~\cite{ACF1} disproved it for $n=3,4$, updating the statement to the current $n\geq 5$, and proved the conjecture for $n=5,6,7,8$, leaving the result open for $n\geq 9$.

\subsection{Symmetric analogues of results in Sections~1-2}
Similar to Example~\ref{E1}, one can prove that the scale by 7 of the fourth symmetric~Gaussian difference ${_q}\Delta_4^s(h,x;f)$ is not symmetric Gaussian, hence ${_q}D_4^sf(x)$ is a counterexample to the Symmetric Gaussian vs. MZ Derivatives Conjecture. In this way, the conjecture is updated to:

\begin{SconjectureG-MZ}[The first update]\label{SCj2} {\rm In orders at least three, each symmetric~MZ differentiation is a scale of a symmetric Gaussian differentiation.}
\end{SconjectureG-MZ}

The following is the symmetric analogue of Proposition~\ref{P8}. Its similar proof is omitted.

\begin{proposition}\label{P8S} {\rm \,(i)} No two distinct symmetric Gaussian differences are scalar multiples of each other.

{\rm (ii)} Two distinct symmetric Gaussian differences are scales of each other if and only~if they have the form
\[\text{${_q}\Delta_n^s$ and either one of ${_{-q}}\Delta_n^s$ and ${_{\pm q^{-1}}}\Delta_n^s$,\quad for $n\geq 3$ and $q\neq 0,\pm 1$.}\]
\end{proposition}

The following is the symmetric analogue of the second update of the Gaussian vs. MZ Derivatives Conjecture.

\begin{SconjectureG-MZ}[The second update]\label{SCj3} {\rm In orders at least three, each symmetric MZ differentiation is equivalent to a symmetric Gaussian differentiation.}
\end{SconjectureG-MZ}

We were unable to provide a counterexample to the first update of the symmetric Gaussian vs. MZ Derivatives Conjecture that is equivalent to a symmetric Gaussian. The reason for this is that no such a counterexample exists. This is highlighted by the following stronger symmetric analogue of Proposition~\ref{P3}.

\begin{proposition}\label{P3S}
The second update of the symmetric Gaussian vs. MZ Derivatives Conjecture is the same as the first update, regardless of exactness.
\end{proposition}

This proposition is a consequence of a stronger symmetric analogue of Theorem~\ref{T1}, the classification of symmetric generalized Riemann derivatives, proved in~\cite{AC2}, which we conveniently restate here in an equivalent form as the following theorem:

\begin{theorem}[\cite{AC2}, Theorem~1.3]\label{T2} \;
Two symmetric generalized Riemann differentiations are equivalent if and only if they are scales of each other.
\end{theorem}

This closes the section and the article. Summarizing, we produced two counterexamples to the Gaussian vs. MZ Derivatives Conjecture, addressed the symmetric analogues of the results,
and updated the conjecture to every MZ derivative is equivalent to a Gaussian. Moving forward, our expectation that the proof of this updated version of the conjecture or of any potential counterexample will require either a new clever idea or a new advanced result in conjunction to the major result of Theorem~\ref{T1} that we used here.

\bibliographystyle{amsplain}

\end{document}